\documentclass[aps,pre,twocolumn,groupedaddress,showpacs,showkeys]{revtex4-1}
\usepackage{graphicx}

\begin{document}

\preprint{Departamento de Física, Universidad Simón Bolívar}

\title{Multidimensional Washboard Ratchet Potentials for Frustrated Two-Dimensional Josephson-Junctions Arrays on square lattices}


\author{Rafael Rangel$^{1,2}$ }
\email[]{rerangel@usb.ve}
\author{ Marcos Negruz$^{1}$}
\affiliation{$^{1}$Departamento de F\'{i}sica, Universidad Sim\'{o}n Bol\'{i}var, A.P. $89000$, Caracas 1080-A, Venezuela}
\affiliation{$^{2}$Laboratorio de F\'{i}sica Te\'{o}rica de S\'{o}lidos(LFTS), Facultad de Ciencias, Universidad Central de Venezuela, A.P. $47586$,  Caracas 1041-A, Venezuela}
\altaffiliation{}




\date{\today}

\begin{abstract}
In this work, we derive an analytical procedure for obtaining  a multidimensional washboard  ratchet potential $U_{\mathbf{f}}$ for two-dimensional
 Josephson junctions array (TDJJA) with an applied magnetic field.
 The magnetic field is given in units of the quantum flux per plaquette or frustration of the form $\mathbf{f=\frac{M}{N}}$. The derivation is done
under the assumption that the checkerboard pattern  ground state or unit cell of a two-dimensional Josephson junctions array (TDJJA) is preserved under current bias.
 The RCSJ with a white noise term models the dynamics for each junction phase in the array. The multidimensional potential is the unique expression of the collective
 effects that emerge from the array in contrast to the single junction. First step in the procedure is to write the equation for the phases for the unit cell by considering
the constraints imposed for the gauge invariant phases due to frustration.
Secondly and key idea of the procedure, is to perform a variable transformation from  the original systems of stochastic equations
to a system of variables where the condition for the equality of mixed second partial is forced via
the Poincar\'{e}'s
 theorem for differentials forms. This leads to a nonlinear matrix equation (equation (9) in the text), where the new coordinates variables $\mathbf{x_{f}}$
are evaluated and where the potential exist.
The transform  matrix also permits the correct transformation of the original white noise terms of each junction to the intensities in the $\mathbf{x_{f}}$ variables.
   The commensurate symmetries of the ground state pinned vortex lattice, leads  to discrete symmetries to the part of the washboard potential  that does not contain  a
tilt due to the external bias current(equation(10) in the text).
In this work  we apply the procedure for the important cases  $\mathbf{f}=\mathbf{\frac{1}{2}},\mathbf{\frac{1}{3}}$. For $\mathbf{f}=\mathbf{\frac{1}{2}}$,
 we show that previously efforts for finding the potential are  restricted, leading to a reduced dimension of the potential.
 Our  potential  reduces to their expression, if one forces their  assumption  which represents an unstable situation. New physics emerge when  currents are applied in the $x$ and $y$ directions,
in particular, we confirm analytically previous
numerical work for $\mathbf{f}=\mathbf{\frac{1}{2}}$, concerning the border of stationary states, a landmark of the potential.
For $\mathbf{f}=\mathbf{\frac{1}{3}}$, we give a generalization of previous  work, in which we include  both the currents in the $x$ and $y$ directions as well the noise terms.
We find the potentials  realize tilted  ratchets analogous to a combustion motor.
\end{abstract}


\maketitle


\section{ Introduction.} The resistive behavior of a two dimensional Josephson array with a given  frustration parameter (the ratio of the perpendicular magnetic field to the flux quantum per plaquette,
 $f=\Phi /\Phi_{0}$, where $\Phi_{0}=hc/2e$),
has been a matter of study
\cite{YongChoiNr1,YongChoiNr2,XSLingNr1,EGranatoNr1,MonTeitelNr1}.
 When the external bias is zero, a mean field approach based quantum interference method can be used to obtain phase diagrams(see \cite{FNoriNr1} and references therein), in which a localization
 without disorder \cite{ReRangelNr1,ReRangelNr2} is exploited. Description of dynamics at any temperature requires knowledge of the origin of dissipation. In superconducting wire networks, near $T_{c}$,
 this is tantamount of using the generalized Ginzburg-Landau equations for each wire element\cite{RTidecks}\cite{FMPeetersNr1}. In Josephson junction arrays the RCSJ model\cite{MTinkham}
 describes each junction. The model contains a tilted washboard potential, that permits to obtain
qualitative  and quantitative understanding of
the dynamics\cite{CLobbNr1,ZondaBelzigNovotny,LuigiLongobardiNr1,SPZhao,JMKinoja,MBorromeo,DevoretEsteveClarke}. Usually these arrays are made such that charging effects  due to small capacitance can be ignored\cite{CLobbNr1}.
In fact, recently interesting questions like switching rates, thermal hopping and retrapping statistics were studied using the analogy of the RCSJ model with the dynamics of
Brownian  particle in a tilted wasboard potential  \cite{ZondaBelzigNovotny}\cite{LuigiLongobardiNr1}. The important question is, if for TDJJA with frustration a similar study can be done, albeit
the high dimensional stotchastic equations have not been formulated.
 On the other hand, TDJJA constitutes natural systems where to realize ratchet behavior \cite{FMPeetersNr2,PeterReimann,HaenggiMarchesoni}, in particular,
 TDJJA ratchets with an asymmetric washboard potential, were fabricated \cite{ShalomPastorizaNr1} in the overdamped regime. Numerical studies \cite{VeronicaMarconi},
  for $f=1/2$, explained some features of the experiments, but qualitative understanding without detailed knowledge of the ratchet landscape  is not possible.
In fact, many researches when pursuing numerical work infer the existence of a multidimensional washboard potential in analogy for
 a single Josephson junction. Others \cite{VShklovskijNr1} used forced uniaxial tilted washboard one-dimensional potentials,
 to model a complex multidimensional physical arrangement \cite{ChallisJack}.\\
\section{TDJJA with Magnetic Fields.}  When  a large $LxL$ network is subject to a uniform  external driving current injected along one
edge of the array and removed at the opposite edge, from energy balance arguments one expects that
the ground state configuration is preserved. Frustration dictates that the net circulation of the vector potential  in a given sense
around a plaquette is given by $2\pi \Phi/\Phi_{0}(mod2\pi)=2\pi(f_{j}-n_{j})$, where $n_{j}$ is an integer which defines the vorticity\cite{MTinkham}.
Neglecting macroscopic screening effects ensures equal frustration for all plaquette $f_{j}=f$. In \cite{FaloBishopLomdahl}, for example the dynamics
is analyzed for $f=\frac{1}{3}$ in the overdamped at $T=0$ and  $T\not= 0$,
 analogous numerical studies are done in \cite{ChungLeeStroud,MonTeitelNr1,VeronicaMarconiNr2,VeronicaMarconiNr3} also in the overdamped regime. Also recently,
interest for the properties of TDJJA near incommensurability have been studied experimentally, in which the frustration value $f=\frac{2}{5}$ appears to be of relevant \cite{YongChoiNr1,YongChoiNr2}.
 Consequently,  we calculate in this work, the multidimensional washboard potential for TDJJA for frustrations $\mathbf{f}=\mathbf{\frac{1}{2}},\mathbf{\frac{1}{3}}$. We manage to find
the multidimensional potential for  the cases $\mathbf{f}=\mathbf{\frac{1}{4}},\mathbf{\frac{1}{5}},\mathbf{\frac{2}{5}}$, however, as
 the complexity and extension of the calculations for rational frustration $\mathbf{f=\frac{M}{N}}$ grow with N, these cases will published in another work.
The implementation of the procedure goes through the following steps:
 1) write the equations for the ground state configuration or basic cell unit, this unit constitutes an $NxN$ array as illustrated in Fig.1  \cite{TeitelJayaprakash,StraleyBarnett}.
 2) identify the variables in order to define a system of stochastic differential equations with diagonal isotropic masses
 3) check that the cross derivatives of the resulting potential are not equal, 4) define a coordinate transformation to a set of new variables, 5) in the newly defined variables,
establish the necessary condition for the existence of a potential by invoking Poincar\'{e}'s theorem for differential forms, 6) we obtain  a non-linear matrix equation (eqt.(9) below),
 whose solution lead us to the potential(eqt.(6) below).
 The multidimensional potentials give the opportunity for studying similar questions posed for a single Josephson junction in theoretical and experimental sense \cite{ZondaBelzigNovotny,LuigiLongobardiNr1,LuigiLongobardiNr2},
 now asked for a TDJJA. We discuss that matter along the work out of the theory below and in the conclusions.\\
\begin{figure}
\includegraphics[width=6cm]{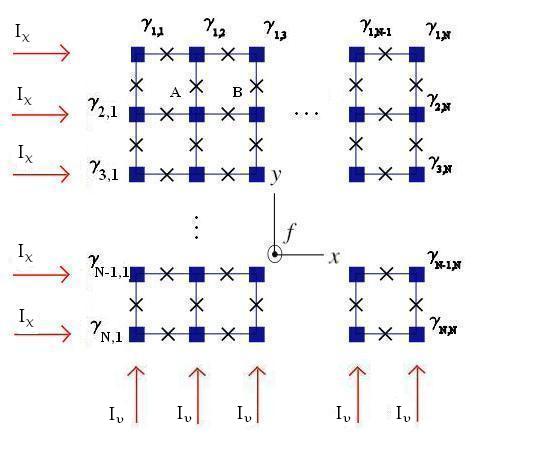}
\caption{NxN Josephson junctions array  unit cell with $f=\frac{M}{N}\Phi_{0}$. The magnetic field is perpendicular to the plane of the sheet and a DC current $I_{\chi}$ is injected in x-direction, and $I_\upsilon$ in the y-direction. The boxes represent superconducting islands, and the crosses the Josephson junctions. }
\label{fig1paper}
\end{figure}
\section{ Equations of Motion.}  For the array in Fig.$(1)$ there are
\cite{TeitelJayaprakash,StraleyBarnett,VeronicaMarconiNr2} $2N(N+1)$  junctions, whose dynamics is given by the RCSJM model\cite{MTinkham}:
\begin{equation}
\beta_{c}\frac{d^{2}\gamma_{i,j}}{d\tau^{2}}+\frac{d\gamma_{i,j}}{d\tau}+\frac{dU(\gamma_{i,j})}{d\gamma_{i,j}}=\sqrt{\frac{2k_{B}T}{E_{J}}}\xi_{i,j}(\tau)
\label{RCSJMEQM}
\end{equation}
where the tilted one-dimensional ``washboard" potential is $U(\gamma_{i,j})=[-\cos(\gamma_{i,j})-\frac{I_{\chi}}{I_{c0}}\gamma_{i,j}]$, and where $\gamma_{i,j}$ is the gauge-invariant phase difference,
$\gamma_{i,j}=\gamma_{i}-\gamma_{j}-(2e/\hbar c)\int_{i}^{j}\mathbf{A}\cdot\mathbf{dl}$,  where $\mathbf{A}$ is the vector potential. $I_{c0}$ is the critical current assumed equal for all junctions.
 $I_{\chi}$ is the current injected in the $x$-direction, and $I_{\upsilon}$ the injected current in the $y$ direction. One has $\tau=(2eI_{c0}R/\hbar)t$ as the dimensionless time, $\beta_{c}=(2eI_{c0}CR^{2}/\hbar)$ is the Stewart--McCumber parameter,
$R$ is the shunt resistance, $C$ is the capacitance, $\hbar$ is the Planck's constant and $e$ is the electron charge.
The stochastic term describes white noise with intensity $\sqrt{\frac{2k_{B}T}{E_{J}}}$,
$\langle \xi_{i,j}(\tau)\rangle=0$,$\langle \xi_{i,j}(\tau)\xi_{i,j}(\tau')\rangle=\frac{2k_{B}T}{E_{J}}\delta(\tau-\tau')$,
 $k_{B}$ is the Boltzmann constant, and
$T$ is the temperature, $E_{J}=(\hbar/2e)I_{c0}$ is the Josephson coupling energy\cite{MTinkham}.
A bookkeeping counting allow us to establish the number of independent equations for the unit cell. First,  periodic boundary conditions imply that the phases in opposite places in Fig.(1) are equal.
 One has $N^{2}$ plaquette, in each of them  the sum of the gauge invariant phases  around each plaquette must be:
$\sum_{plaquete}\gamma_{i}=2\pi f(mod2\pi)=2\pi(f-n_{j})$, when going around a contour, vorticity  is given by:
$\sum_{contour}\gamma_{i}=2\pi\sum_{enclosed-cells}(f-n_{j})$, ground state symmetries constrain further the number of phases. First, the circulation around the perimeter of the unit cell,
should be zero, in order to avoid size scale dependent energy terms. Secondly, the circulation around a contour formed from plaquette in any column or any row should also be zero.
For $f=\frac{1}{N}\Phi_{0}$, the zero circulation can always be achieved in any contour around a column or row by putting
has one vortex in a selected plaquette ( circulation$=-2(\pi \frac{N-1}{N})$, and zero vortices in the other ($(N-1)$ plaquette ( circulation$=2(\pi \frac{1}{N})$ en each).
The particular pattern configuration  of the minimal energy, i.e., the $n_{j}$, were found in \cite{TeitelJayaprakash} for the frustrations we are interested.
 One has also, that only different phases are the ones in the first column and the first row.
The internal phases are just suitable combinations of these phases, i.e., there are $2N$ different phases. We calculate them from the current conservation
 in the $(N-1)$ internal nodes in the second row ( symbols A,B...) and current conservation in the $x$ and $y$ directions. We have also $(N+1)$ current equations for $2N$ unknowns.
We use then the $(N-1)$ relations from the plaquette circulation to eliminate the remaining  $(N-1)$  phases. One finally has $(N+1)$ current equations for the  $(N+1)$ independent variables $y_{i}$,
 they make the functions $g_{j}(y_{1},\ldots,y_{\small{N+1}})$ (primes are derivatives with respect the dimensionless time, $\beta_{c}$ is the Stewart-McCumber parameter):
\begin{equation}
\beta_{c}y_{j}''+y_{j}'+g_{j}(y_{1},\ldots,y_{N+1})=0; \hspace{0.5cm} j=1,\ldots,N+1
\label{firstequationmotion}
\end{equation}

\begin{widetext}
\begin{equation}
g_{j}(y_{1},\ldots,y_{N+1})=\frac{1}{2}\sum_{k=1}^{2N}\omega_{jk}\sin\Phi_{k}(y_{i})-\frac{I_{\chi}}{2I_{c}}\delta_{j1}-\frac{I_{\upsilon}}{2I_{c}}\delta_{j2}
+\sqrt{\frac{2k_{B}T}{E_{J}}}\left [\sum_{m=1}^{2}\delta_{jm}\sum_{l=1}^{N}\xi_{l}(\tau)+\sum_{m=3}^{N+1}\delta_{jm}\sum_{l=1}^{4}\xi_{l}(\tau)\right ]
\label{two}
\end{equation}
\end{widetext}
$y_{1}$ is sum of the phases in the left side, and $y_{2}$  is the sum of the phases in the upper side in fig.(1).
The rest of the variables $y_{i},\small{i=3,..,N+1}$,  are chosen, in consistence with the form of eqt.(\ref{firstequationmotion}).
Furthermore, $\delta_{jm}$ is the Kroenecker delta. The original  phases $\Phi_{k}, \small{k=1,..2N}$ are now  function of the $(N+1)$ new variables
 $y_{i}$ and the $2Nx(N+1)$ matrix $\omega_{jk}$,  with entries  $-1,0,1$ gives  the presence of the functions $\sin\Phi_{k}(y_{i})$. In eqt.(\ref{two}),
for the variable $y_{1}$ there are $2N$ independent noise terms, similar for the variable $y_{2}$. For the each of the remaining variables $y_{i}, i=3,N+1$,
 there are four independent noise terms. One finds  $\partial g_{j}/\partial y_{i}\neq\partial g_{i}/\partial y_{j}$, i.e., $g_{j}$ are not the derivatives of a potential\cite{HarleyFlanders}.
 One searches for a new set of variables $\mathbf x$, and a transformation  is carried out through an $(N+1)x(N+1)$ matrix $\mathbf D$, $\mathbf{ x=Dy}$, and
corresponding inverse transformation $\mathbf{y=D^{-1}x}$. Multiplying equation (\ref{firstequationmotion}) with $D_{ij}$, one obtains a new system of stochastic differential equations:

\begin{widetext}
\begin{equation}
\beta_{c}\underbrace{\sum_{j=1}^{N+1}D_{ij}y_{j}''}_{x_{i}''}+\underbrace{\sum_{j=1}^{N+1}D_{ij}y_{j}'}_{x_{i}'}+\underbrace{\sum_{j=1}^{N+1}D_{ij}g_{j}(y_{1},y_{2},\ldots,y_{N+1})}_{f_{i}(x_{1},x_{2},\ldots,x_{N+1})}=0
\label{cua}
\end{equation}
\end{widetext}

\begin{widetext}
\begin{equation}
f_{i}=\sum_{k=1}^{2N}a_{ik}\sin\Phi_{k}(x_{j})-\frac{I_{\chi}}{2I_{c}}D_{i1}-\frac{I_{\upsilon}}{2I_{c}}D_{i2}+\sqrt{\frac{2k_{B}T}{E_{J}}}\left [\sum_{m=1}^{2}D_{im}\sum_{l=1}^{N}\xi_{l}(\tau)+\sum_{m=3}^{N+1}D_{im}\sum_{l=1}^{4}\xi_{l}(\tau)\right ]
\label{definitionfi}
\end{equation}
\end{widetext}

where $a_{ik}=\frac{1}{2}\sum_{j=1}^{N+1}D_{ij}\omega_{jk}$. One needs to find  a defining equation for the matrix $\mathbf D$, such that in the new variables the cross derivatives are equal.
This is achieved in the next section.
\section{Construction of The Potential.} We define 1--form  by $F=f_{i}dx_{i}$, $i=1,2,\ldots,N+1$, and force the condition $dF=0$, i.e, the 1-form is closed. We invoke the Poincar\'e's theorem and look for a  0--form $U$,
 such that $dU=F$, implying $d(dU)=0$, i.e, the 1--form is exact\cite{HarleyFlanders}, i.e., we obtain the potential. Define $\Omega=\sqrt{\frac{2k_{B}T}{E_{J}}}$. The form of  $f_{i}$,
suggest the following  Ansatz for the 0--form:
\begin{widetext}
\begin{equation}
U(x_{i},\Omega,I_{\chi},I_{\upsilon})=-\sum_{k=1}^{2N}\cos\Phi_{k}(x_{i})-\frac{I_{\chi}}{2I_{c}}\sum_{i=1}^{N+1}D_{i1}x_{i}-\frac{I_{\upsilon}}{2I_{c}}\sum_{i=1}^{N+1}D_{i2}x_{i}
+\Omega\left [\sum_{i=1}^{N+1}\sum_{m=1}^{2}D_{im}x_{i}\sum_{l=1}^{N}\xi_{l}(\tau)+\sum_{i=1}^{N+1}\sum_{m=3}^{N+1}D_{im}x_{i}\sum_{l=1}^{4}\xi_{l}(\tau)\right ]
\label{och}
\end{equation}
\end{widetext}
\begin{widetext}
\begin{equation}
dU=\sum_{i=1}^{N+1}\left [\sum_{k=1}^{2N}\frac{\partial\Phi_{k}}{\partial x_{i}}\sin\Phi_{k}-\frac{I_{\chi}}{2I_{c}}D_{i1}-\frac{I_{\upsilon}}{2I_{c}}D_{i2}+\Omega(\sum_{m=1}^{2}D_{im}\sum_{l=1}^{2N}\xi_{l}(\tau)+\sum_{m=3}^{N+1}D_{im}\sum_{l=1}^{4}\xi_{l}(\tau))\right ]dx_{i}=F;
\label{nuv}
\end{equation}
\end{widetext}

therefore, from eqt.(\ref{definitionfi}), the necessary condition for $dU=F$ is given by: $a_{ik}=\frac{\partial\Phi_{k}}{\partial x_{i}}$,

\begin{equation}
d(dU)=\sum_{i=1}^{N+1}\sum_{m=1}^{N+1}\frac{\partial}{\partial x^{m}}(dU)dx^{m}\wedge dx^{i}=0,
\label{ele}
\end{equation}

A necessary and sufficient condition for $d(dU)=dF=0$, is the equality of mixed partial derivatives:
$\frac{\partial f_{i}}{\partial x_{m}}=\frac{\partial f_{m}}{\partial x_{i}}; \hspace{0.5cm} i\neq m$, we use the last relation in obtaining  $\mathbf D$, note that last condition
 can be written as $a_{ik}\frac{\partial \Phi_{k}}{\partial x_{m}}=a_{mk}\frac{\partial \Phi_{k}}{\partial x_{i}}$,
by noting that $\Phi_{\partial x,k}=\Phi_{\partial y,k}D^{-1}$, one obtains:

\begin{equation}
\frac{1}{2}D^{T}D=\Phi_{\partial y,k}\omega^{T}(\omega\omega^{T})^{-1}
\label{thi}
\end{equation}
where $\Phi_{\partial y,k}=\frac{\mathbf{\Phi_{k}}}{\mathbf{\partial y}}$ represents the matrix of the derivatives of phases with respect to the variables $y_{j}$.
After solving eqt.(8) for $D$, we read out the potential (eqt.(6)) and the equations of motions(eqts.$(4-5)$). Equivalently, the functions $f_{i}(x_{j})$,
which define the equations of motion (equation$(4)$),  can be found by differentiating the potential:
\begin{equation}
f_{i}(x_{1},...,x_{N+1})=\frac{\partial U}{\partial x_{i}}; \hspace{0.5cm} i=1,...,N+1
\end{equation}
 The potential  has in general a
periodic part with period $\vec{a}$ and a linear tilt, i.e.,
\begin{eqnarray}
U(\vec{r})=U_{\mathbf{0}}(\vec{r})-\vec{g}\cdot\vec{r}\nonumber\\
 \vec{g}=(g_{1},g_{2},0,..0)\nonumber\\
 U_{\mathbf{0}}(\vec{r})=U_{\mathbf{0}}(\vec{r}+a_{j}\hat{\mathbf{r}}_{j})=U_{\mathbf{0}}(\vec{r}+\vec{a})\nonumber\\
\label{tilt}
\end{eqnarray}
 Therefore one has:
\begin{widetext}
\begin{equation}
U(x_{1}+\delta x_{1},x_{2}+\delta x_{2},..,x_{N+1}+\delta_{N+1})=U(x_{1},..x_{N+1})+\delta U_{1}+\delta U_{2}
\end{equation}
\end{widetext}
 where $\delta x_{i}$ is the period in direction $i$, and $\delta U_{1,2}$
is the increment of the potential due to the applied currents.
In a regime when the noise terms and dissipation can be neglected, one obtains a Hamiltonian system with
\begin{equation}
H= \sum_{i=1}^{N+1}\frac{\xi^{2}_{i}}{2\beta_{c}}+U(x_{1},..x_{N+1})
\end{equation}
The fist term is the kinetic energy, $\frac{\xi^{2}_{i}}{2\beta_{c}}=\frac{\beta_{c}}{2}\dot{x}^{2}_{i}$(the dot represents the time derivative).
 The Stewart and McCumber parameter can be written in the form:
\begin{equation}
\frac{\beta_{c}}{2}=\frac{4E_{j}}{E_{C}}\frac{R^{2}}{R^{2}_{Q}}
\end{equation}
 where $E_{J}$ is the Josephson coupling energy already  previously defined, $E_{C}=\frac{e^{2}}{2C}$ is the charging energy, and $R_{Q}=\frac{2\hbar}{2e}$, is the quantum resistance.
Equation $(13)$ is the generalization of the one junction case (see eqt.$(4)$ in \cite{CLobbNr1},\cite{MTinkham}) as the potential is a multidimensional one. Quantization of eqt.$(12)$ is an standard task,
the relevance of which is given for the case when $E_{C}\geq E_{J}$ and $R\geq R_{Q}$. On the other hand, the existence of minima of the potential warrant  the stability of the quantum system.
We calculate the border of stability as a function of the applied currents for the fist example we discuss in the next section.
 Without noise but maintaining the dissipation terms, one can analyze the flow properties of the associated first order system:
\begin{equation}
\beta_{c} X^{'}= \Xi
\end{equation}
\begin{equation}
\Xi^{'}= -\Xi/\beta_{c}-\nabla U(x_{1},..,x_{N+1})
\end{equation}

 where $X=(x_{1},..,x_{N+1})$, $\Xi=(\xi_{1},..,\xi_{N+1})$\cite{SWiggins}, in this case the dynamical system is phase space contracting \cite{LichtenbergLieberman}.
 For $\beta_{c}=0$, one has the overdamped limit,
\begin{equation}
 X^{'}= -\nabla U(x_{1},..,x_{N+1})
\end{equation}
a gradient system with at most fixed points
 \cite{GuckenheimerHolmes}. Proper interplay of nonlinearities and noise in the systems we derive below, happens in the underdamped  regime.\\
\section{Potential for $\mathbf{f=1/2}$.}  Consider the ground state $2\times2$ superlattice unit cell for $f=1/2$ in Fig.$(2)$( see Fig.$(1)$ in \cite{TeitelJayaprakash}):
\begin{figure}
\includegraphics[width=6cm]{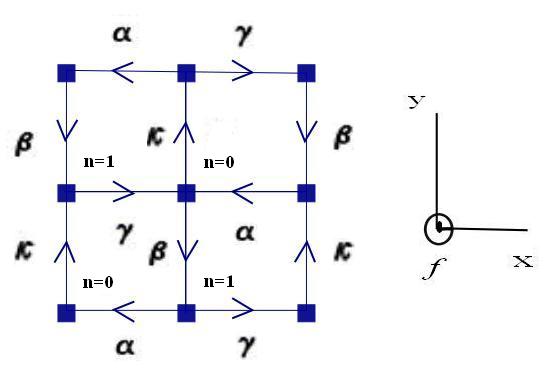}
\caption{A 2x2 Josephson junctions array under a magnetic field $f=\frac{1}{2}\Phi_{0}$}
\label{fig3paper}
\end{figure}
 Following section III, one first obtains three equations from current conservation in the $x$, $y$ direction and from the central node:
\begin{widetext}
\begin{eqnarray}
\beta_{c}(\gamma''-\alpha'')&+&(\gamma'-\alpha')+(\sin\gamma-\sin\alpha)=I_\chi/I_{c}\nonumber\\
\beta_{c}(\kappa''-\beta'')&+&(\kappa'-\beta')+(\sin\kappa-\sin\beta)=I_\upsilon/I_{c}\\
\beta_{c}(\gamma''+\alpha''-\beta''-\kappa'')&+&(\gamma'+\alpha'-\beta'-\kappa')+(\sin\gamma+\sin\alpha-\sin\beta-\sin\kappa)=0\nonumber
\label{cat}
\end{eqnarray}
\end{widetext}
Then, secondly, the quantization condition  for the n=0 plaquette is given by: $\beta+\kappa+\alpha+\gamma=\pi$, it allows to eliminate the $(\alpha + \gamma )$ variables.
 We arrive to the choice of variables: $\mathbf{y}=(y_{1},y_{2},y_{3})$ as:
\begin{eqnarray}
y_{1}=\frac{1}{2}(\gamma-\alpha)\nonumber\\
y_{2}=\frac{1}{2}(\kappa-\beta)\\
y_{3}=-(\beta+\kappa)\nonumber
\label{sev}
\end{eqnarray}

We follow  eqts.(3-5) for $N=2$.  The $\Phi_{k}$ variables
 are:
$\Phi_{k}(y_{i})=(\Phi_{1}=\alpha,\Phi_{2}=\beta,\Phi_{3}=\gamma,\Phi_{4}=\kappa)$, and the matrix  $\omega$, which can be read from eqt.$(18)$, is given by:

\[ \omega = \left( \begin{array}{cccc}
-1 & 0 & 1 & 0 \\
0 & -1 & 0 & 1 \\
1 & -1 & 1 & -1
\end{array} \right)\]

one has:
\[ \Phi_{\partial y,k} = \left( \begin{array}{cccc}
-1 & 0 & 1 & 0 \\
0 & -1 & 0 & 1 \\
1/2 & -1/2 & 1/2 & -1/2
\end{array} \right)\]

First one proves that $\partial g_{j}/\partial y_{i}\neq \partial g_{i}/\partial y_{j}$, and proceed to calculate the right hand side of eqt.(9), one obtains:

\[ D^{T}D = \left( \begin{array}{ccc}
2 & 0 & 0 \\
0 & 2 & 0 \\
0 & 0 & 1
\end{array} \right)\]

This implies,

\[ D = \left( \begin{array}{ccc}
\sqrt{2} & 0 & 0 \\
0 & \sqrt{2} & 0 \\
0 & 0 & 1
\end{array} \right)\]

and its inverse,

\[ D^{-1} = \left( \begin{array}{ccc}
1/\sqrt{2} & 0 & 0 \\
0 & 1/\sqrt{2} & 0 \\
0 & 0 & 1
\end{array} \right)\]

The new variables, $\mathbf{ x=Dy}$ are:

\begin{equation}
x_{1}=\sqrt{2}y_{1} ;x_{2}=\sqrt{2}y_{2} ;x_{3}=y_{3}; \hspace{0.5cm}
\end{equation}

We define $x=x_{1},y=x_{2},z=x_{3}$ and obtain with eqt.$(5)$, the potential $U$,
\begin{widetext}
\begin{equation}
U(x,y,z,\Omega,I_{\chi},I_{\upsilon})=-\sum_{k=1}^{4}\cos\Phi_{k}(x,y,z)+\frac{\sqrt{2}}{2}(\Omega\sum_{l=1}^{2}\xi_{l}(\tau)-\frac{I_{\chi}}{I_{c}})x+
\frac{\sqrt{2}}{2}(\Omega\sum_{l=3}^{4}\xi_{l}(\tau)-\frac{I_{\upsilon}}{I_{c}})y+\frac{1}{2}\Omega\sum_{l=5}^{8}\xi_{l}(\tau)z
\end{equation}
\end{widetext}
The phases $\Phi_{k}(x,y,z)$  in eqt.($18)$ are given by:
\begin{eqnarray}
\Phi_{1}=\alpha=\frac{1}{2}(\pi-\sqrt{2}x+z)\nonumber\\
 \Phi_{2}=\beta=-\frac{1}{2}(\sqrt{2}y+z)\nonumber\\
\Phi_{3}=\gamma=\frac{1}{2}(\pi+\sqrt{2}x+z)\nonumber\\
\Phi_{4}=\kappa=\frac{1}{2}(\sqrt{2}y-z)\nonumber\\
\end{eqnarray}

One has from eqt.(\ref{tilt}): $U(\vec{r})=U_{\mathbf{0}}(\vec{r})-\vec{g}\cdot\vec{r}$, with
$\vec{r}=x_{i}\hat{r}_{i}$, $i=1,2,3$, $|\hat{r}_{i}|=1$, $\vec{g}=\frac{1}{\sqrt{2}}\frac{I_{\chi}}{i_{c}}\hat{r}_{1}+\frac{1}{\sqrt{2}}\frac{I_{\upsilon}}{i_{c}}\hat{r}_{2}$,
and period $\vec{a}=a_{i}\hat{r}_{i}=4\pi(\frac{1}{\sqrt{2}},\frac{1}{\sqrt{2}},1)$.
\begin{figure}
\includegraphics[width=6cm]{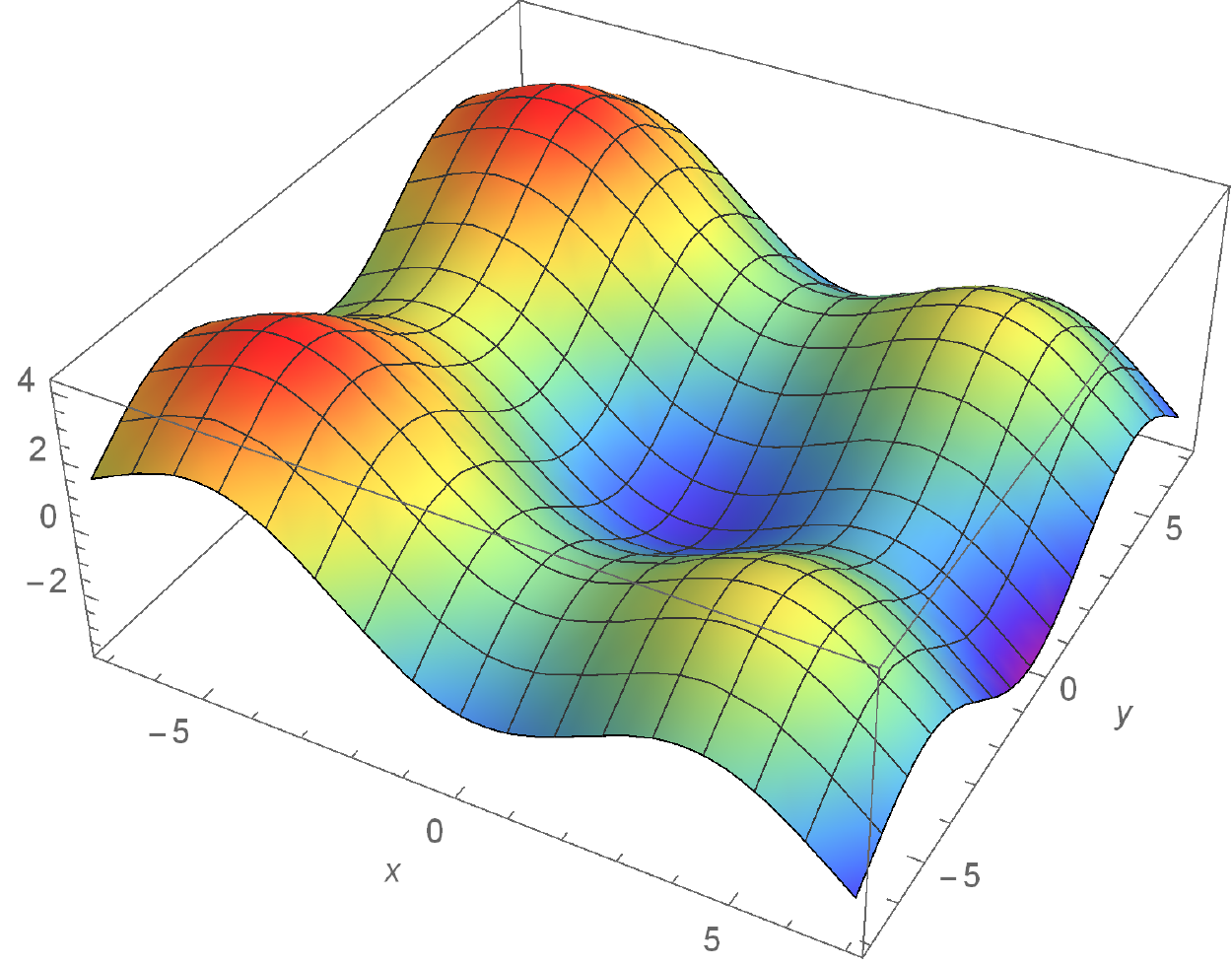}
\caption{$U(x,y,\pi/2)$ from eqt.(13),$\Omega=0$, $I_{\chi}=0.292893$,  $I_{\upsilon}=0.0$}
\end{figure}
\begin{figure}
\includegraphics[width=6cm]{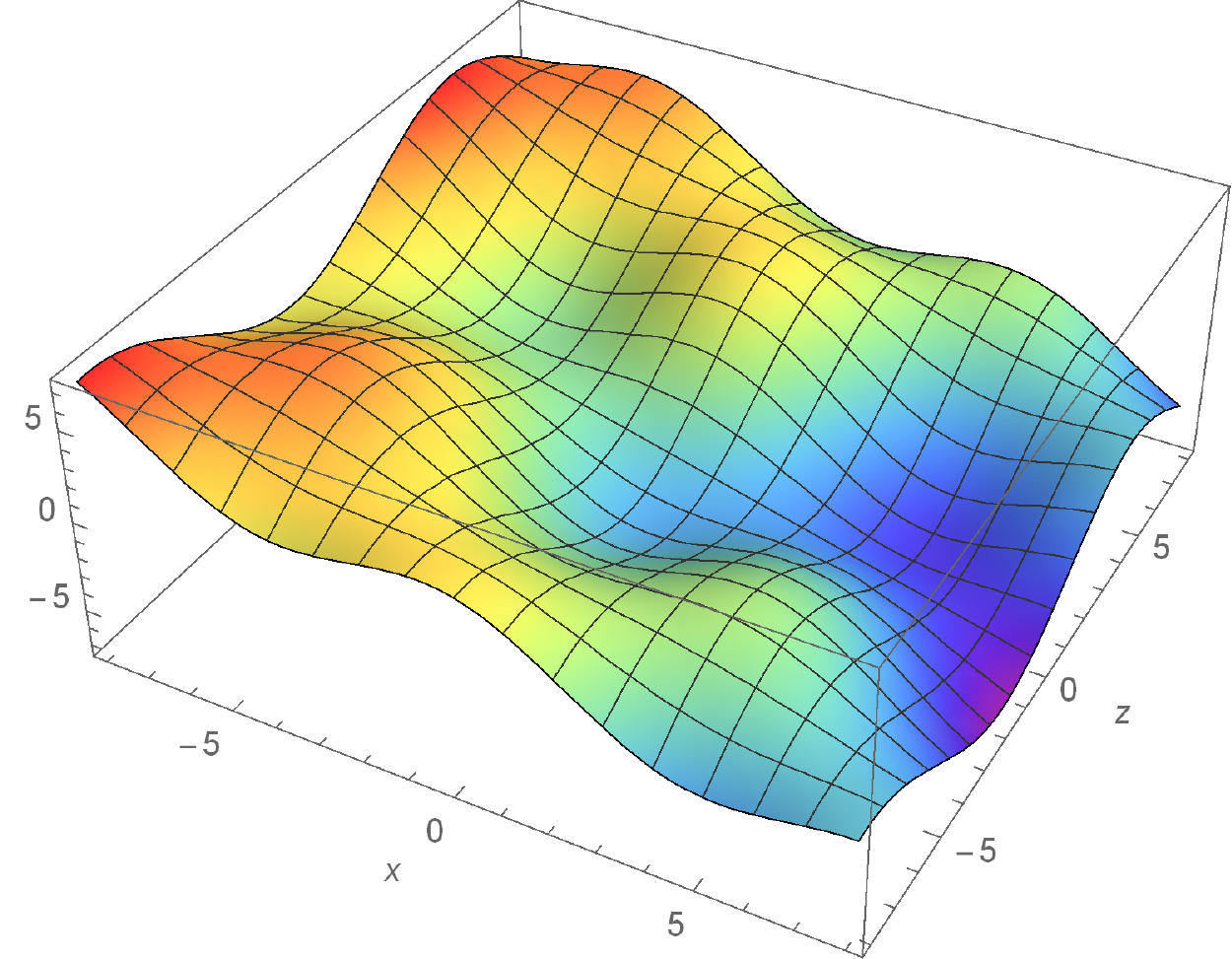}
\caption{$U(x,0,z)$ from eqt.(13),$\Omega=0$, $I_{\chi}=0.8485$,  $I_{\upsilon}=0.0$}
\end{figure}
\begin{figure}
\includegraphics[width=6cm]{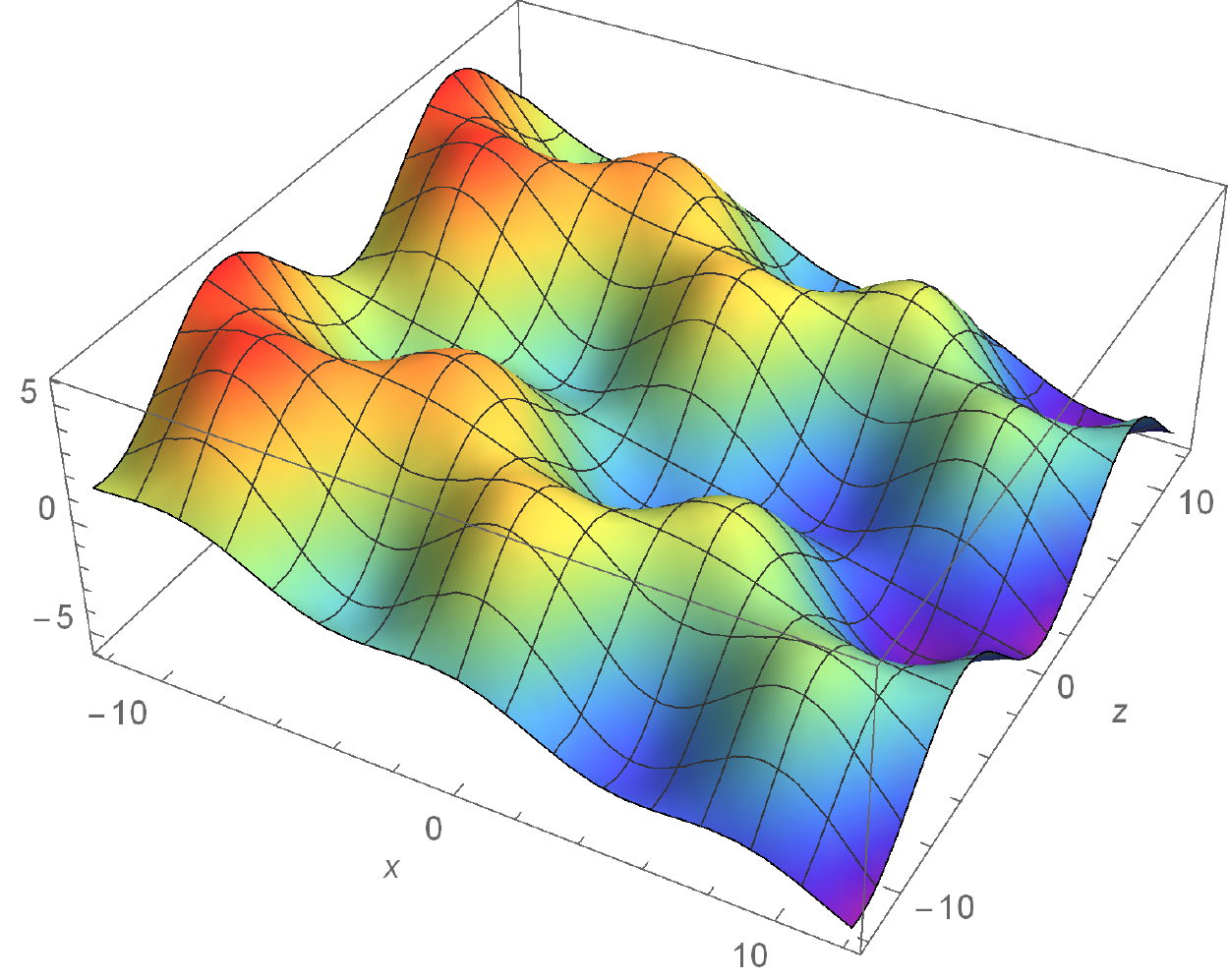}
\caption{$U(x,0,z)$ from eqt.(13),$\Omega=0$, $I_{\chi}=0.4142$,  $I_{\upsilon}=0.0$}
\end{figure}
\begin{figure}
\includegraphics[width=6cm]{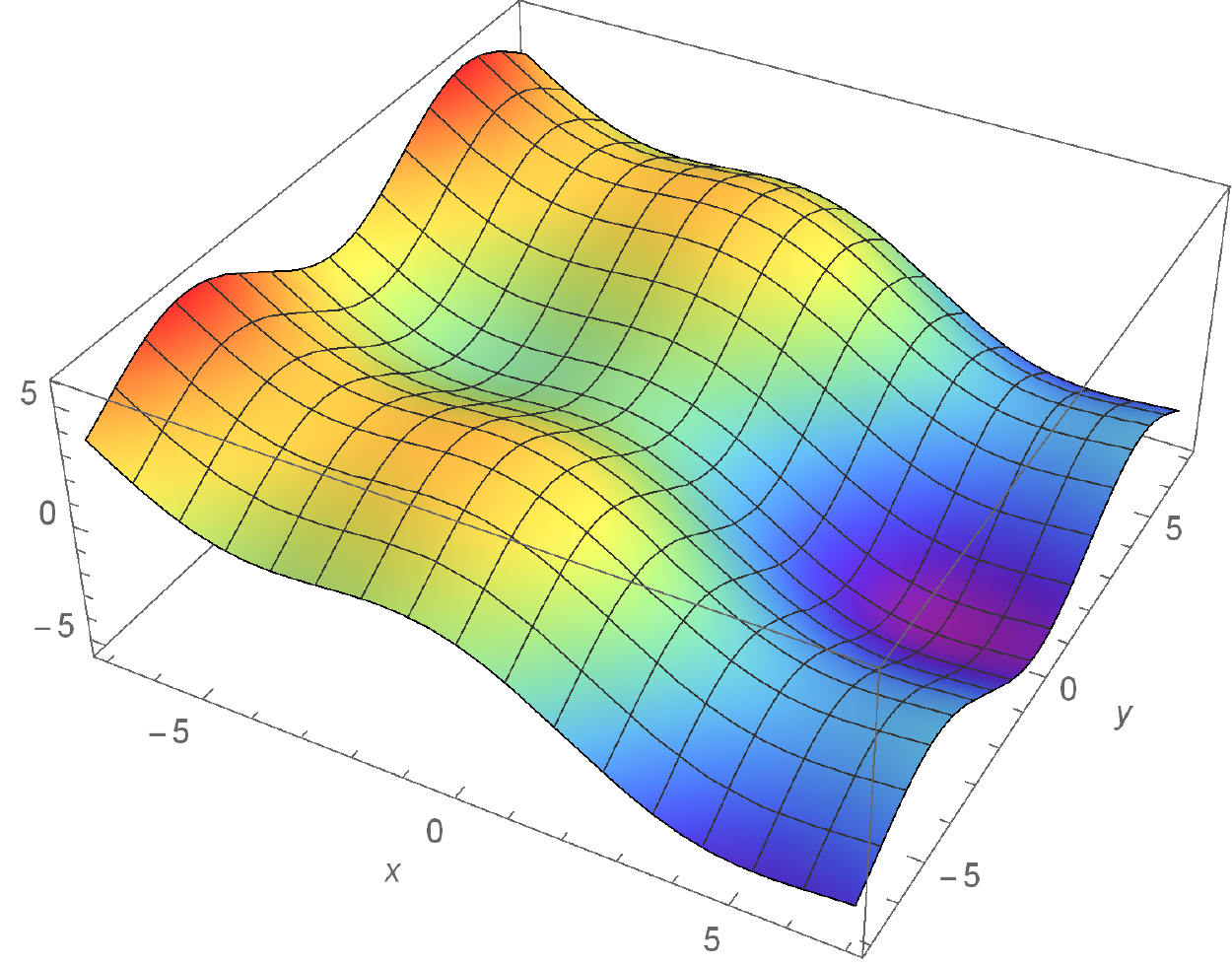}
\caption{$U(x,y,-\pi/2)$ from eqt.(13),$\Omega=0$, $I_{\chi}=0.828429$,  $I_{\upsilon}=0.0$}
\end{figure}
\begin{figure}
\includegraphics[width=6cm]{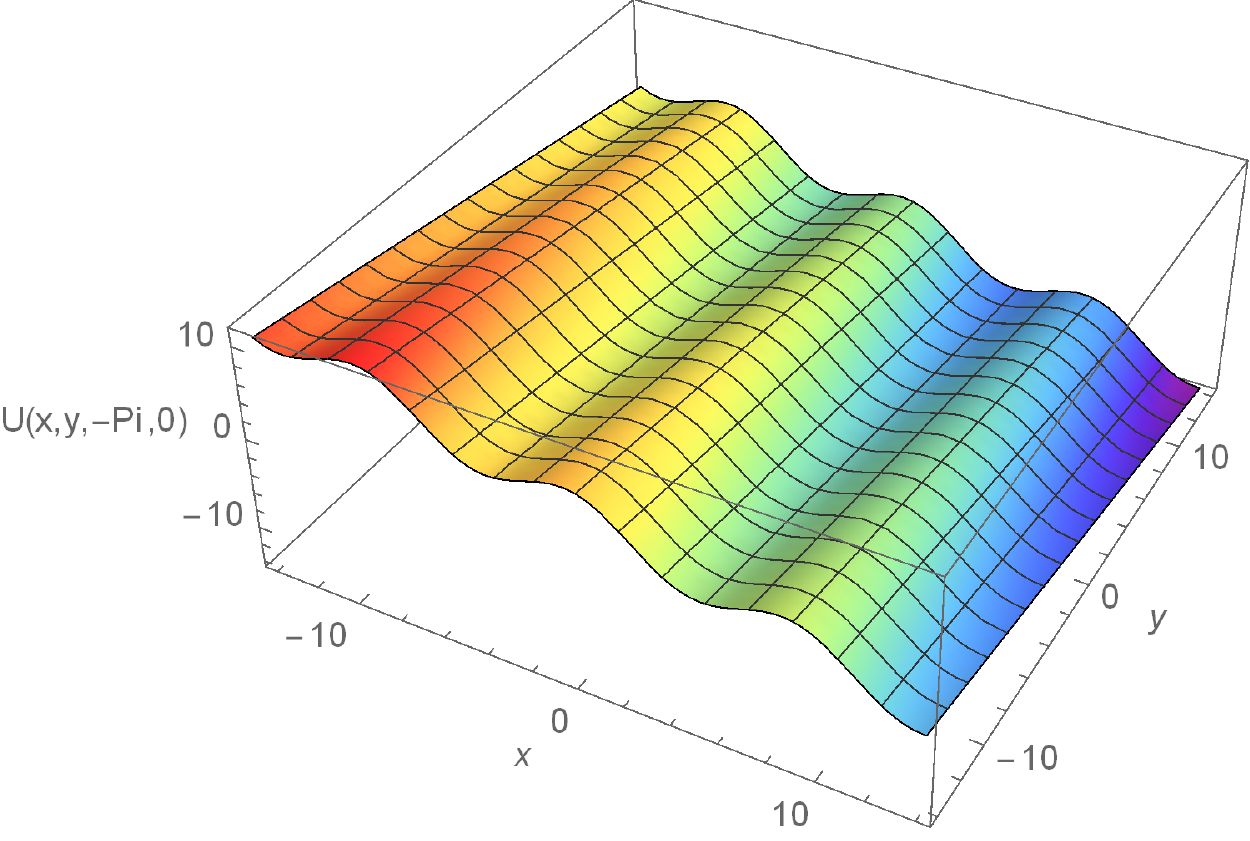}
\caption{$U(x,y,-\pi)$ from eqt.(13),$\Omega=0$, $I_{\chi}=0.849942$,  $I_{\upsilon}=0.3999$}
\end{figure}
\begin{figure}
\includegraphics[width=6cm]{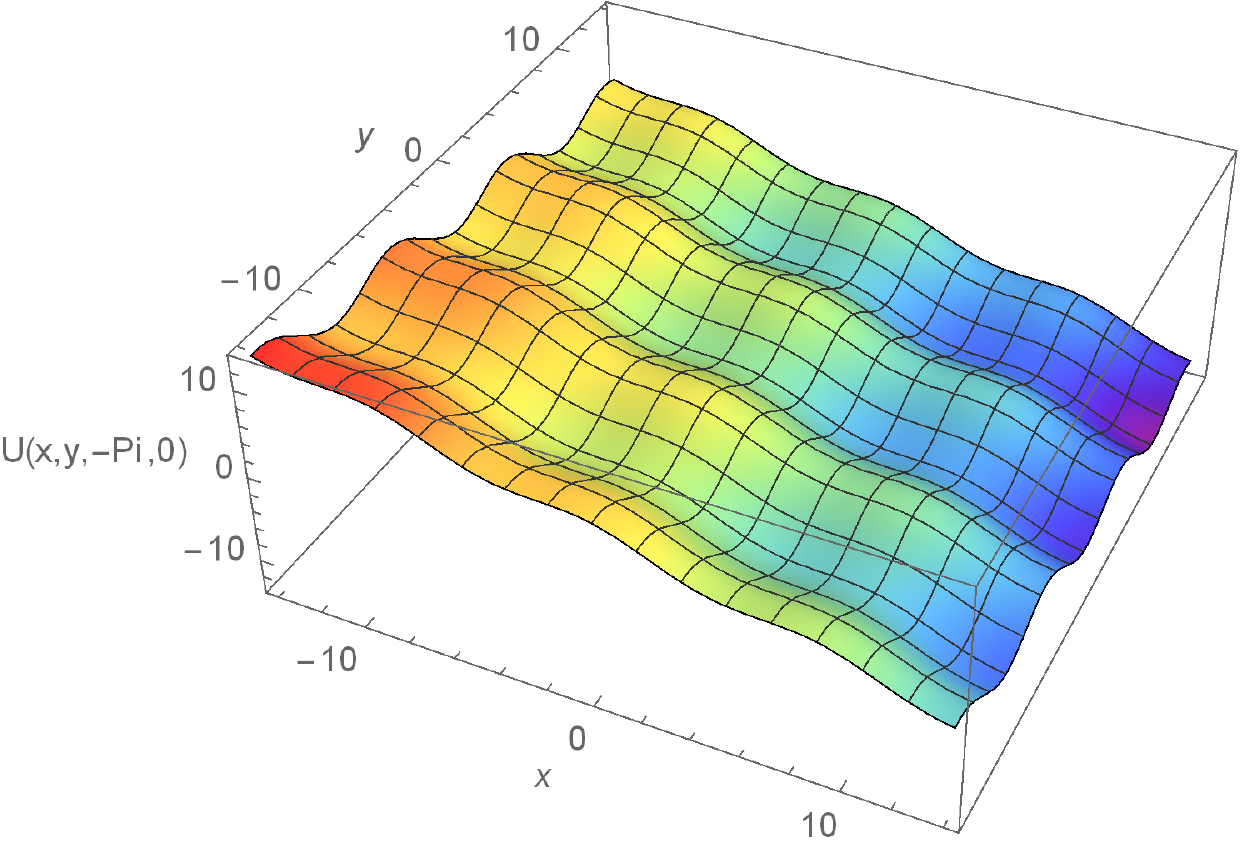}
\caption{$U(x,y,-\pi/4)$ from eqt.(13),$\Omega=0$, $I_{\chi}=0.8499$,  $I_{\upsilon}=0.3999$}
\end{figure}
Figures$(3-8)$ show the projection of the potential $U(x,y,z)$ for some specific values of the variables y or z.

With the knowledge of $\mathbf D$ the corresponding equations of motion (eqts.$(4)$ and$(5)$) can be straightforwardly written, or alternatively, the functions $f_{i}(x,y,x)$ can be
obtained using eqt.$(10)$:
\begin{eqnarray}
\beta_{c}x''&+&x'-\frac{2}{\sqrt{2}}\sin(\frac{z}{2})\sin(\frac{x}{\sqrt{2}})= \frac{1}{\sqrt{2}}(\frac{I_{\chi}}{I_{c0}}+\Omega\sum_{l=1}^{2}\xi_{l})\nonumber\\
\beta_{c}y''&+&y'+\frac{2}{\sqrt{2}}\cos(\frac{z}{2})\sin(\frac{y}{\sqrt{2}})= \frac{1}{\sqrt{2}}(\frac{I_{\upsilon}}{I_{c0}}+\Omega\sum_{l=1}^{2}\xi_{l})\nonumber\\
\beta_{c}z''&+&z'+\cos(\frac{z}{2})\cos(\frac{x}{\sqrt{2}})+\sin(\frac{z}{2})\cos(\frac{y}{\sqrt{2}})=\frac{\Omega}{2}\sum_{l=1}^{4}\xi_{l}\nonumber\\
\label{cat}
\end{eqnarray}
For $\Omega=0$, there is a stationary time-independent regime. Stable solutions in this regime correspond to local minima of the potential. In this case, for given $I_{\chi}, I_{\upsilon}$,
one manipulates eqt.$(23)$, and obtain a relation which defines the $x$ variable:
\begin{equation}
\frac{\sin(\sqrt{2}x)}{\sqrt{\frac{1}{2}(1+\sqrt{1-(\frac{I_{\upsilon}}{I_{\chi}})^{2}\sin^{2}(\sqrt{2}x)})
+\cos^{2}(\frac{x}{\sqrt{2}})}}=I_{\chi},
\end{equation}
$z(I_{\chi})$ is obtained from the relation $\tan(z/2)=-\cos(x/\sqrt{2})/\cos(y/\sqrt{2}))$ and $y(I_{\chi})$ from $\sin(x/\sqrt{2}))/\sin(y/\sqrt{2}))=(\frac{I_{\upsilon}}{I_{\chi}})^{2}$.
 For $I_{\upsilon}=0$, there is a critical current, at the value of which, the local minima and maxima merge, i.e.,
at the critical current the stable and unstable fixed points  coalesce into one. This critical fixed point  can be obtained by maximizing eqt.(24)
 with $I_{\upsilon}=0$, one obtains: $x^{crit}=1/\sqrt{2}\cos^{-1}{(2\sqrt{2}-3)}$, $y^{crit}=0$, $z^{crit}=-2\sin^{-1}{\sqrt{(2-\sqrt{2})/2}}$,
 and $I^{crit}_{\chi}=\sqrt{2\sqrt{2}}\sqrt{3\sqrt{2}-4}$ per cell, in agreement with the calculation done in \cite{Ben}.
Identical results are obtained, as one should expect by symmetry, if one put $I_{\chi}=0$ and $I_{\upsilon}\neq 0$.
For $I_{\chi¨}<I^{crit}_{\chi}$, the local minima of the potential occurs always  at $y=0$, an example is shown in fig.$(3)$. Also local maxima of the potential happens at $y=0$.
For $I_{\chi}>I^{crit}_{\chi}$, however, the potential has no local minimum,
which implies the non-existence of fixed points. Recall that $y=0$ implies $\kappa=\beta$, if one uses this constraint in the potential (eqt.(21)), one obtains another system with
one dimension reduced. This system exists in the projection of the potential to the line $y=0$. Next we show this statement.
Suppose now that the terms representing stochastic forces are irrelevant ($\Omega\approx 0$) and take the imposed transversal current $I_{\upsilon}$ to be  zero (see fig.(1)). If one assume
one could pin the variable $y$ to the value zero at any current (see fig. (2)),
in this case, the potential reduces to:

\begin{equation}
U(x,0,z)=-\cos\alpha-2\cos\beta-\cos\gamma-\frac{I_{\chi}}{\sqrt{2}I_{c}}x
\end{equation}
where $\alpha,\beta$ and $\gamma$ are given by

\begin{eqnarray}
\Phi_{1}=\alpha&=&\frac{1}{2}(\pi-\sqrt{2}x_{1}+x_{3})\label{fifty}\\
\Phi_{2}=\beta&=&-\frac{1}{2}x_{3}\\
\Phi_{3}=\gamma&=&\frac{1}{2}(\pi+\sqrt{2}x_{1}+x_{3})\label{fiftytwo}
\end{eqnarray}

This potential can be rewritten using trigonometric identities and the sum and the difference of the equations (\ref{fifty}) and (\ref{fiftytwo})

\begin{widetext}
\begin{equation}
U(\alpha,\gamma)=-2\cos(\frac{\alpha+\gamma}{2})\cos(\frac{\alpha-\gamma}{2})-2\sin(\frac{\alpha+\gamma}{2})+\frac{I}{I_{c}}(\frac{\alpha-\gamma}{2})
\end{equation}
\end{widetext}

Now we introduce the scaled sum and difference variables $\xi\equiv (\alpha+\gamma)/\sqrt{2}$ and $\eta\equiv(\alpha-\gamma)/2$,  this allows  to write the potential as

\begin{equation}
U=2[-\cos(\frac{\xi}{\sqrt{2}})\cos(\eta)-\sin(\frac{\xi}{\sqrt{2}})+\frac{I}{2I_{c}}\eta]=2U(\xi,\eta)
\end{equation}
One obtains the system:
\begin{eqnarray}
\beta_{c}\eta''&+&\eta'+\frac{\partial U(\xi,\eta)}{\partial \eta}\nonumber\\
\beta_{c}\phi''&+&\phi'+\frac{\partial U(\xi,\eta)}{\partial \phi}\nonumber\\
\end{eqnarray}

This system was used for simulations in \cite{Oct}. Forcing  $y=0$ for $I_{\chi}>I^{crit}_{\chi}$ however,
 is incompatible with the phase flow properties  of eqt.(10)\cite{SWiggins}, since for currents greater than the critical, the line $y=0$ is neither local nor global attracting. Instead,
 there will limit cycles and  fluctuations in the variable $y$, i.e., voltage fluctuations in the $y$-direction \cite{VeronicaMarconiNr4}.  In the overdamped regime $\beta_{c}=0$,
 for  $I_{\upsilon}\neq 0$, and $I_{\chi}\neq 0$, Fisher et tal., carried out
 numerical calculations for $\Omega=0$\cite{FisherStroudJanin}, as equations of motion they used eqt.$(18)$ for the gauge invaraint phases $\Phi_{k}$, in the overdamped regime (eqt.$(17)$).
They found a regime with voltage zero, that they called a pinned regime. The border of this regime can be obtained  analytically from the calculation of the maximum permitted value of $I_{\upsilon}$
 for a given $I_{\chi}$. Again, by maximizing eqt.$(24)$ one finds a polynomial  equation of degree six for the unknown $\cos (\sqrt{2} x)$. For a given  $ R \equiv I_{\upsilon}/I_{\chi}$,
the solution of the polynomial equation allows to calculate the maximal value of $I_{\chi}$, i.e., the border of stability of the pinned phase. For  the special value $R =1$,
 the algebraic equation reduces to one of degree four, one finds the solution $x=y=\pi /2\sqrt{2}$,
\begin{figure}
\includegraphics[width=6cm]{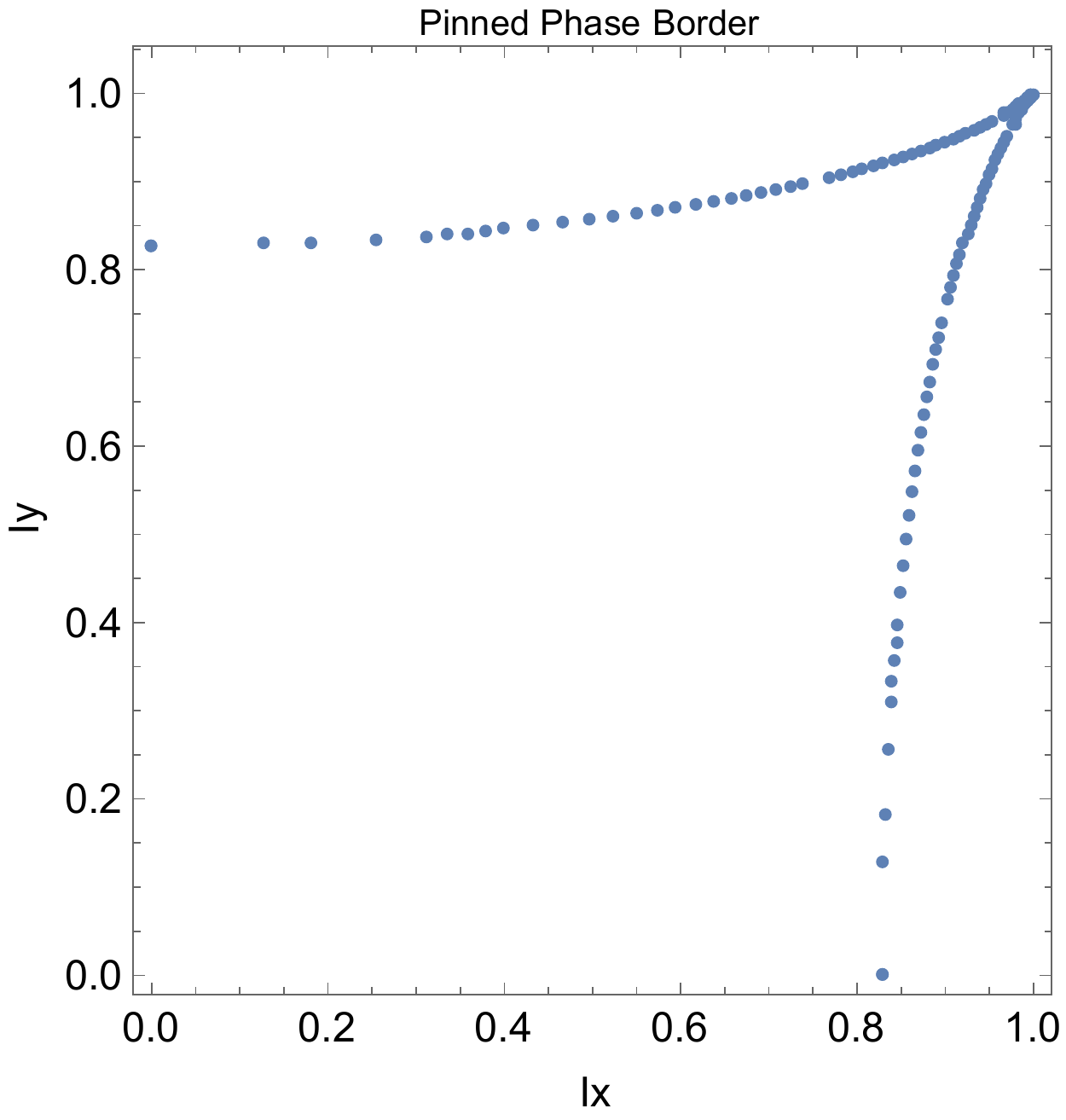}
\caption{The pinned phase obtained from eqt.(24).}
\end{figure}
 and $z=-\pi /2$, which implies $I^{max}_{\chi}=1$ per cell. When $R =0$ one has the case discussed above. Therefore, there is an island of stability between
 $I_{\chi}^{cri}$ and $I_{\chi}^{max}$. In order to get an analytical equation for the border of stability, one finds first from eqt.$(24)$ a polynomial equation of degree four for
 $y\equiv \cos (\sqrt{2} x)$,

\begin{widetext}
\begin{equation}
y^{4}+ I^{2}_{\chi}y^{3}+[2(I^{2}_{\chi}-1)+\frac{I^{4}_{\chi}}{4}(1-R^2)]y^{2}+[(I^{2}_{\chi}-1)I^{2}_{\chi}]y+[(I^{2}_{\chi}-1)^{2}-\frac{I^{4}_{\chi}}{4}(1-R^2)]=0
\end{equation}
\end{widetext}
One factorizes the $y=-1$ root, and obtains:
\begin{widetext}
\begin{equation}
y^{3}+[(I^{2}_{\chi}-1)]y^{2}+[(I^{2}_{\chi}-1)+\frac{I^{4}_{\chi}}{4}(1-R^2)]y+[(1-2I^{2}_{\chi})+I^{4}_{\chi}(1-\frac{(1-R^{2})}{4}]=0
\end{equation}
\end{widetext}
This equation gives for given $I_{\chi}$ and $R$ the corresponding value of $x,y,z$. Third,
one writes the discriminant $D(I_{\chi},R)$, of this third degree equation and look for its change in sign, i.e., one solves $D(I_{\chi},R)=0$ which for given $I_{\chi}$,
 is a polynomial equation of third degree
for $R_{min}^{2}(I_{\chi})$, the value of which defines the pinned phase border (the value of the branch of one of the  real roots for $\cos (\sqrt{2} x)$ that evolves from $R=1 $
 to the value $R_{min}$,
defined as the value of $R$ where it transform into a complex root).
 Fig.$(9)$ shows the pinned phase. This phase is a landmark property of the potential independent of $\beta_{c}$.
Fig.(6) in \cite{FisherStroudJanin} shows a numerical calculation of this  exact analytical result, the difference of  factor two in the axis is because we use the notion of critical current per cell.
Beyond the pinned border of the stationary regime, there are no local minima of the potential and only time dependent solutions exist. For finite temperatures,  numerical simulation seems to be
the only way to study this regime, however, qualitative understanding can be obtained from the potential. The kind of questions one can ask were put in \cite{VeronicaMarconiNr3}. The authors made
numerical simulations  of large arrays with frustration $f=1/25$ in the overdamped limit. Their phase diagram Temperature versus applied current, (their fig.(1) show various phases). The pinned phase,
with not voltage corresponds to  the stationary regime shown in fig.$(5)$, it destabilizes for sufficiently large value of $\Omega$, transforming it in a phase with a finite time average voltage.
The mechanism behind is  similar to well known case of a single Josephson junction \cite{AmbegoakarHalperin,LuigiLongobardiNr1}. Only that the barrier height $\Delta U$, has to be calculated from eqt.(10) with $\Omega=0$
and the criterium that the scape rate turns significant when $ \Delta U \approx \Omega$ \cite{FerrandoSpadaciniTommei}. In this way, tilting  the potential asymmetrically, i.e., by applying currents below or above the $I_{\chi}=I_{\upsilon}$ line,
it is clear that one direction destabilizes first, and then for larger $\Omega$, the other direction.  This is only a qualitative picture, and a quantitative analysis  needs the full nonlinear dynamics
and particular properties of the potential in order to understand the final state after escaping. Also proper use of a multidimensional Wiener process \cite{KDebrabantARoessler} is required.
On the other hand, the so called transverse depinning \cite{VeronicaMarconiNr3}, viewed form our theory constitutes the ratchet effect similar to previous case. When $I_{\chi}\ge I_{crit}$ and one turns
 $I_{\upsilon} \simeq\epsilon $ on, one begins to tilt the potential in the $y$ direction. There are channels  around $z=n*\pi$, where the potential permits the particle of mass $\beta_{c}$ to slide almost free down, whereas
for example $z=\pi/4$ it is halted by a relatively  big potential barrier(see figures (7-8)), then at at sufficiently big value of $\Omega$,
the particle begins to slide down the direction $y$, accompanied with a voltage.
 What we have at hand is the analogous of a molecular motor \cite{Mar,KellerBustamante}. The numerical study of these scenarios, also eventually the combination of a constant current and time dependent
periodic current is a matter of research. This last scenario has been treated numerically  in \cite{VeronicaMarconiNr2} for $f=\frac{1}{2}$   at zero temperature ($\Omega =0$) in the overdamped regime ($\beta_{c} =0$).\\
\section{Potential for  $\mathbf{f=1/3}$.} In the appendix, we apply the general method (equations $5$ and $9$), to that case and obtain the potential.
Deeper analytical and numerical work of that case is challenging future task.\\
\section{ Conclusions.} In this letter, we have developed a general method to find a potential for current biassed  TDJJA with frustration. We analyzed in some detail
the frustration value $f=1/2$ for which new analytical
results are found, one important result is the analytical calculation of the pinned  phase, as it has a landmark character deriving
 from the potential. In the overdamped regime, a rocking ratchet effect was found in \cite{VeronicaMarconi}, where an asymmetrical potential was engineered. For our potential, we conjecture
 that inertial effects ($\beta_{c}\ne 0$) can produce a dynamical ratchet. Also we expect the current reversal phenomenon to exist in our system \cite{VincentMayerKurths}, in fact our systems are
prominent examples of ratchets in which inertia, dissipation and noise combines together with high dimensionality and chaos\cite{RangelNr3} similar to the theory of molecular motors
 \cite{Mar,KellerBustamante}\cite{PeterReimann}\cite{HaenggiMarchesoni}.
 Potentials for symmetric values of f around $f=1/2$, for example $f=1/4$ and $f=3/4$ \cite{VeronicaMarconi}, appear interesting and can also be found with the method. At temperatures
$T\simeq 0$, when dissipation can be neglected, i.e., neglecting quasi-particles degrees of freedom \cite{RangelVerrilliMarin}, ours systems are Hamiltonian ones as explained (see equation $(13)$).
In this case \cite{GSchoen,LuigiLongobardiNr2}, the potentials derived here still can be used in conjunction with quantum noise \cite{Gardinir1}. If the charging energy due to a small capacitance is comparable to the
Josephson coupling energy $E_{J}$ \cite{CLobbNr1} , the problem turns a quantum mechanical one\cite{MTinkham}\cite{LuigiLongobardiNr2}. With the aim of studying superconducting to non-superconducting transition
 for TDJJA with current bias but
no magnetic field, Porter and Stroud \cite{PorterStroud} writes a Hamiltonian  similar to our eqt.$(13)$, in which the kinetic terms are charging energies, which in the units used in \cite{MTinkham}, is the energy
of a particle of mass $C(\hbar/2e)^{2}$, and the potential is the sum of the one-dimensional washboard  potentials of the junctions.
 The search is for local minima of an unknown
 multidimensional potential as a signal of the superconducting phase. Due to the lack of dissipation and fluctuations the only possible transition is from superconducting to non-supeconduncting states.
 This task for the TDJJA with frustration can be studied as we have shown with the help of
 the multidimensional potentials derived here. Finally, we mention the potentially interesting question posed in \cite{GSchoenNr2,WolfgangBelzig} concerning the nature of the fluctuations generated by a
single Josephson junction. We believe the study of the same question for our multidimensional system is a relevant issue.
\begin{acknowledgments}
M.N. acknowledges support from the DID, USB and the ACFMN
 during part of the course of this work.
\end{acknowledgments}

\appendix
\section{Potential for $\mathbf{f=1/3}$.} 
\begin{figure}
  \includegraphics[width=6cm]{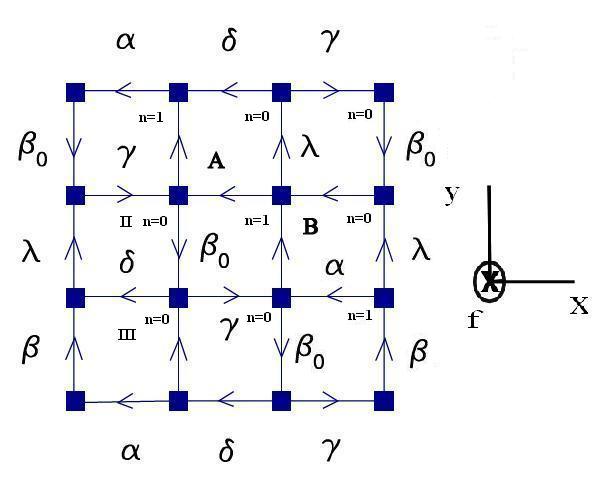}\\
  \caption{A 3x3 Josephson junctions array under a magnetic field $f=\frac{1}{3}\Phi_{0}$}\label{fig3neg}
\end{figure}

Consider the ground state 3x3 superlattice unit cell for $f=1/3$ as in figure (\ref{fig3neg})\cite{TeitelJayaprakash,RangelNr3}.\\

Like in the previous case,  we derive the equations of motion for this arrangement. First we write the flux quantization conditions from the plaquette labeled
II and III in figure \ref{fig3neg}, second we write the conservation of charge conditions at nodes A and B. Then we write the equations of currents in $x$ and $y$ directions.
Later we introduce new variables to obtain an isotropic mass tensor. From the resulting system of equations we read the matrices  $\omega$ and $\Phi_{\partial y,k}$,
form the right hand side of eqt.$(9)$ and find the matrix $\mathbf D$. Then we write the potential (eqt.$(6)$ and the equations of motion (eqts.$(4-5)$.
We have the following flux quantization conditions from the plaquette labeled II and III:

\begin{eqnarray}
\lambda+\gamma+\beta_{0}+\delta=\frac{2\pi}{3} \hspace{0.5 cm} \label{c1}\\
\alpha+\beta-\delta-\lambda=\frac{2\pi}{3} \hspace{0.5 cm} \label{c2}
\end{eqnarray}

From current conservation of charge, we obtain from nodes A and B:

\begin{widetext}
\begin{eqnarray}
\beta_{c}(\gamma''+\alpha''-\beta''-\beta_{0}'')+(\gamma'+\alpha'-\beta'-\beta_{0}')+(\sin\gamma+\sin\alpha-\sin\beta-\sin\beta_{0})=0\label{c3}\\
\beta_{c}(\beta''-\alpha''+\delta''-\lambda'')+(\beta'-\alpha'+\delta'-\lambda')+(\sin\beta-\sin\alpha+\sin\delta-\sin\lambda)=0\label{c4}
\end{eqnarray}
\end{widetext}

We impose the condition for the currents in the $x$ $y$ direction
\begin{equation}
\beta_{c}(\gamma''-\alpha''-\delta'')+(\gamma'-\alpha'-\delta')+(\sin\gamma-\sin\alpha-\sin\delta)=I_{\chi}/I_{c}
\end{equation}
\begin{equation}
\beta_{c}(\beta''-\beta_{0}''+\lambda'')+(\beta'-\beta_{0}'+\lambda')+(\sin\beta-\sin\beta_{0}+\sin\lambda)=I_{\upsilon}/I_{c}
\end{equation}
The conditions (\ref{c1}) and (\ref{c2}) allow us to rewrite the equations (\ref{c3}) and (\ref{c4}) respectively

\begin{widetext}
\begin{eqnarray}
2\beta_{c}(\gamma''+\alpha'')+2(\gamma'+\alpha')+(\sin\gamma+\sin\alpha-\sin\beta-\sin\beta_{0})=0\\
2\beta_{c}(\beta''-\lambda'')+2(\beta'-\lambda')+(\sin\beta-\sin\alpha+\sin\delta-\sin\lambda)=0
\end{eqnarray}
\end{widetext}

We choose the variables $\mathbf{y}=(y_{1},y_{2},y_{3},y_{4})=(x,y,z,u)$. The scaling in the definitions of the new variables is necessary to obtain a diagonal and isotropic mass tensor

\begin{eqnarray}
x&=&1/2(\gamma-\alpha-\delta)\nonumber\\
y&=&1/2(\beta-\beta_{0}+\lambda)\nonumber\\
z&=&1/2(\beta-\alpha+\delta-\lambda)=\beta-\lambda-\pi/3\label{c9}\\
u&=&1/2(\gamma+\alpha-\beta-\beta_{0})=\gamma+\alpha-2\pi/3\nonumber
\end{eqnarray}

Using the equations (\ref{c1}), (\ref{c2}) and (\ref{c9}), the gauge invariant phases $\Phi=(\alpha,\beta,\beta_{0},\gamma,\lambda,\delta)$ can be obtained

\begin{eqnarray}
\alpha&=&1/3(u-z-2x+\pi)\nonumber\\
\beta&=&1/3(2y-u+z+\pi)\nonumber\\
\beta_{0}&=&1/3(-2y-2u-z+\pi)\nonumber\\
\gamma&=&1/3(2u+z+2x+\pi)\nonumber\\
\lambda&=&1/3(2y-u-2z)\\
\delta&=&1/3(u+2z-2x)\nonumber
\end{eqnarray}

So we can write this system in compact manner

\begin{equation}
\beta_{c}y_{j}''+y_{j}'+g_{j}(y_{1},y_{2},y_{3},y_{4})=0\label{c11}
\end{equation}

with $j=1,2,3,4$; $\mathbf{y}=(y_{1}=x,y_{2}=y,y_{3}=z,y_{4}=u)$,\newline$\Phi_{\mathbf k}=(\Phi_{1}=\alpha,\Phi_{2}=\beta,\Phi_{3}=\beta_{0},\Phi_{4}=\gamma,\Phi_{5}=\lambda,\Phi_{6}=\delta)$
and

\begin{equation}
g_{j}(y_{1},..,y_{4})=\frac{1}{2}\sum_{k=1}^{6}\omega_{jk}\sin\Phi_{k}(y_{1},..,y_{4})\\
-\frac{I_{\chi}}{2I_{c}}\delta_{j1}-\frac{I_{\upsilon}}{2I_{c}}\delta_{j2}
\end{equation}

$\delta_{ji}$ is the kronecker delta and,


\[ \omega = \left( \begin{array}{cccccc}
-1 & 0 & 0 & 1 & 0 & -1 \\
0 & 1 & -1 & 0 & 1 & 0  \\
-1 & 1 & 0 & 0 & -1 & 1 \\
1 & -1 & -1 & 1 & 0 & 0 \\
\end{array} \right)\]

\[ \Phi_{\partial y,k} = \left( \begin{array}{cccccc}
-2/3 & 0 & 0 & 2/3 & 0 & -2/3 \\
0 & 2/3 & -2/3 & 0 & 2/3 & 0 \\
-1/3 & 1/3 & -1/3 & 1/3 & -2/3 & 2/3 \\
1/3 & -1/3 & -2/3 & 2/3 & -1/3 & 1/3
\end{array} \right)\]

Using equation $(9)$ in the main text one obtain matrix $D$:

\[ D = \left( \begin{array}{cccc}
2/\sqrt{3} & 0 & 0 & 0 \\
0 & 2/\sqrt{3} & 0 & 0 \\
0 & 0 & 1 & 1 \\
0 & 0 & 1/\sqrt{3} & -1/\sqrt{3}\\
\end{array} \right)\]

\[ D^{-1} = \left( \begin{array}{cccc}
\sqrt{3}/2 & 0 & 0 & 0  \\
0 & \sqrt{3}/2 & 0 & 0 \\
0 & 0 & 1/2 & \sqrt{3}/2 \\
0 & 0 & 1/2 & -\sqrt{3}/2  \\
\end{array} \right)\]

From the inverse transformation $D^{-1}$ one  writes:

\begin{eqnarray}
y_{1}&=&\frac{\sqrt{3}}{2}x_{1}\nonumber\\
y_{2}&=&\frac{\sqrt{3}}{2}x_{2}\nonumber\\
y_{3}&=&\frac{1}{2}x_{3}+\frac{\sqrt{3}}{2}x_{4}\nonumber\\
y_{4}&=&\frac{1}{2}x_{3}-\frac{\sqrt{3}}{2}x_{4}\nonumber\\
\end{eqnarray}

With eqt.$(6)$ one reads the potential:

\begin{widetext}
\begin{eqnarray}
U&=&-\sum_{k=1}^{6}\cos\Phi_{k}(x_{1},x_{2},x_{3},x_{4})+\frac{1}{\sqrt{3}}(\Omega\sum_{l=1}^{3}\xi_{l}(\tau)-\frac{I_{\chi}}{I_{c}})x_{1}+\frac{1}{\sqrt{3}}(\Omega\sum_{l=4}^{6}\xi_{l}(\tau)
-\frac{I_{\upsilon}}{I_{c}})x_{2}+ \frac{1}{2}\Omega\sum_{l=7}^{14}\xi_{l}(\tau)x_{3}\nonumber\\
&+&\frac{1}{2\sqrt{3}}\Omega[\sum_{l=7}^{10}\xi_{l}(\tau)- \sum_{l=11}^{14}\xi_{l}(\tau)]x_{4}
\end{eqnarray}
\end{widetext}

Also the equations of motion in the new variables $\mathbf x$ can be written just by reading the corresponding elements of the matrix D as eqt.$(5)$ dictates, or from the potential
using eqt.$(10)$.
The gauge invariant phases as a function of the new variables $\mathbf{ x=Dy}$, in eqt.$(5)$  and in eqt.$(A13)$ are:

\begin{eqnarray}
\Phi_{1}=\alpha&=&1/3(-\sqrt{3}x_{1}-\sqrt{3}x_{4}+\pi)\nonumber\\
\Phi_{2}=\beta&=&1/3(\sqrt{3}x_{2}+\sqrt{3}x_{4}+\pi)\nonumber\\
\Phi_{3}=\beta_{0}&=&1/3(-\sqrt{3}x_{2}-\frac{3}{2}x_{3}+\frac{\sqrt{3}}{2}x_{4}+\pi)\nonumber\\
\Phi_{4}=\gamma&=&1/3(\frac{3}{2}x_{3}-\frac{\sqrt{3}}{2}x_{4}+\sqrt{3}x_{1}+\pi)\nonumber\\
\Phi_{5}=\lambda&=&1/3(\sqrt{3}x_{2}-\frac{3}{2}x_{3}-\frac{\sqrt{3}}{2}x_{4})\nonumber\\
\Phi_{6}=\delta&=&1/3(\frac{3}{2}x_{3}+\frac{\sqrt{3}}{2}x_{4}-\sqrt{3}x_{1})\nonumber\\
\end{eqnarray}
This potential can be written according to equation (\ref{tilt}) as

\begin{equation}
U(\vec{r})=U_{\mathbf{0}}(\vec{r})-\vec{g}\cdot\vec{r}
\end{equation}
with $\vec{r}=x_{i}\hat{r}_{i}$, $i=1,\ldots,4$, $|\hat{r}_{i}|=1$,
 $\vec{g}=\frac{1}{\sqrt{3}}\frac{I_{\chi}}{i_{c}}\hat{r}_{1}+\frac{1}{\sqrt{3}}\frac{I_{\upsilon}}{i_{c}}\hat{r}_{2}$
and period
\begin{equation}
\vec{a}=a_{i}\hat{r}_{i}=\frac{2}{\sqrt{3}}\pi(3,3,2,6)
\end{equation}

The functions $f_{i}(x_{1},x_{2},x_{3},x_{4})$ from their definition in equation $(4)$, are  given by:
\begin{widetext}
\begin{eqnarray}
f_{1}&=&\frac{2}{\sqrt{3}}g_{1}=\frac{1}{\sqrt{3}}(\sin\gamma-\sin\alpha-\sin\delta+ \Omega\sum_{l=1}^{3}\xi_{l}(\tau)-I_{\chi}/I_{c})\nonumber\\
f_{2}&=&\frac{2}{\sqrt{3}}g_{2}=\frac{1}{\sqrt{3}}(\sin\beta-\sin\beta_{0}+\sin\lambda+ \Omega\sum_{l=4}^{6}\xi_{l}(\tau)-I_{\upsilon}/I_{c})\nonumber\\
f_{3}&=&g_{3}+g_{4}=\frac{1}{2}(\sin\delta-\sin\lambda+\sin\gamma-\sin\beta_{0}+ \Omega\sum_{l=7}^{14}\xi_{l}(\tau))\nonumber\\
f_{4}&=&\frac{1}{\sqrt{3}}(g_{3}-g_{4})=\frac{1}{2\sqrt{3}}(2\sin\beta-2\sin\alpha+\sin\delta-\sin\lambda-\sin\gamma+\sin\beta_{0}+ \Omega\sum_{l=7}^{10}\xi_{l}(\tau)- \Omega\sum_{l=11}^{14}\xi_{l}(\tau))\nonumber\\
\end{eqnarray}
\end{widetext}
In this way, the equations of motion for the $x_{i}$ variables (eqt.$(4)$) are completed. The border of the pinned phase, in analogy with the $f=1/2$ case,
 and other questions of the
dynamics  will be accomplished  in another work.

\end{document}